\documentstyle[12pt,amscd,amssymb]{article}

\def\su#1{{\sp{(#1)}}} 
\def\tens{\mathop{\otimes}}

\def\<{{\langle}}
\def\>{{\rangle}}

\def\note#1{{}}

\def\note#1{}

\def\ba{{\bar{a}}}

\def\Bar{{\rm Bar}}
\def\ob{{\rm Obs}}
\def\binom#1#2{\left(\begin{array}{c}#1\\ #2\end{array}\right)}

\def\Label{\label}
\begin{document}

\newtheorem{proposition}{Proposition}[section]
\newtheorem{lemma}[proposition]{Lemma}
\newtheorem{corollary}[proposition]{Corollary}
\newtheorem{theorem}[proposition]{Theorem}
\newtheorem{example}[proposition]{Example}

\newcommand{\Section}{\setcounter{definition}{0}\section}
\newcounter{c}
\renewcommand{\[}{\setcounter{c}{1}$$}
\newcommand{\etyk}[1]{\vspace{-7.4mm}$$\begin{equation}\Label{#1}
\addtocounter{c}{1}}
\renewcommand{\]}{\ifnum \value{c}=1 $$\else \end{equation}\fi}
\renewcommand\refname{\begin{center}\vspace{-\baselineskip}\normalsize REFERENCES
\end{center}}

~{\ }\qquad\hskip 3.8in \hskip -0.62802pt \phantom{~}\vspace{1in} 
\begin{center} 
{\large\bf DEFORMATION OF ALGEBRA FACTORISATIONS} 
\vspace{\baselineskip}\\ 
{\sc Tomasz Brzezi\'nski}\footnote{EPSRC Advanced Research Fellow.}\\
{ Department of Mathematics,  University of Wales Swansea,\\
 Singleton Park, Swansea SA2 8PP, U.K.}\\  \& \\
   Department of Theoretical Physics, University of \L\'od\'z,\\ 
Pomorska 149/153, 90--236 \L\'od\'z, Poland. \\
{E-mail}: {\tt T.Brzezinski@swansea.ac.uk}
\end{center}
\vspace{3\baselineskip}

\baselineskip = 17pt 
\begin{abstract}
A general deformation theory of algebras which factorise into two subalgebras is
studied. It is shown that the classification of deformations is related
to the cohomology of a certain double complex reminiscent of the
Gerstenhaber-Schack complex of a bialgebra.
\end{abstract}
\section{Introduction}
An algebra factorisation or a twisted tensor product is a unital algebra $X$ over a field $k$ 
together with two (unital) subalgebras
$B$, $A$ such that the map $B\otimes A\to X$ given by multiplication is
an isomorphism. In what follows we identify $X$ with
$B\otimes A$ as a $(B,A)$-bimodule via this isomorphism and hence the
algebra structure on $X$ can be viewed as a twisting of the usual tensor
product algebra. The algebra $X$ consists of elements
of the form $x=\sum_ib_i\otimes a_i =\sum_ib_ia_i$, where $b_i\in B$ and
$a_i \in A$. An
algebra factorisation is denoted by $X(B,A)$. Algebra factorisations
appear frequently in algebra and number theory.  Examples include the
 tensor product and the braided tensor
product algebras. Also
the quaternions can be
viewed as an algebra factorisation over the real numbers built on two 
copies of the complex numbers (cf.\ \cite{BrzMaj:fac}\cite{CaeIon:fac}). 
Of physical interest are
algebra factorisations obtained by the quantisation of phase spaces. For
example, the
Heisenberg algebra is a factorisation built on the
algebras generated by the momentum and  position operators. 

The theory of deformations of a twisted tensor product between
{\em undeformed} algebras $B$ and $A$ was introduced in \cite{CapSch:twi},
where it has been shown that the Heisenberg algebra is such a 
deformation of the usual tensor product.  
The aim of this note is to give a cohomological interpretation of 
a general theory of
deformations of algebra factorisations which allows for deformation of 
all the algebras
$A$, $B$ and $X$ entering the factorisation.  

In \cite{Tam:coe}\cite{Maj:fou}\cite{CapSch:twi} it is shown that factorisations are in one-to-one
correspondence with linear maps $\Psi:A\otimes B\to B\otimes A$ such
that
\begin{equation}
\Psi(\mu_A\tens B)=(B\tens\mu_A)(\Psi \otimes
A)(A\otimes \Psi),
\quad \Psi(1_A\tens B)=B\tens 1_A,
\label{diag.A}
\end{equation}
\begin{equation}
\Psi(A\tens\mu_B)=(\mu_A\tens B)(B\otimes
\Psi)(\Psi\otimes B),
\quad \Psi(A\tens 1_B)=1_B\tens A.
\label{diag.B}
\end{equation}
Here and below the following notation is used. For an algebra $A$, 
the identity map on $A$ is denoted by $A$,  the
unit of $A$, viewed either as an element of $A$ or as a map $k\to A$
is denoted by $1_A$ and the product map is denoted by $\mu_A$.
For a given factorisation $X(B,A)$ the map $\Psi$ is given by $ab =
(\mu_X\circ\Psi)(a\otimes b)$, and the equations (\ref{diag.A}),
(\ref{diag.B}) simply express 
the associativity conditions: $(aa')b = a(a'b)$, $a(bb')=(ab)b'$ and the
unit condition $1_Ab = b$, $a1_B=a$, for all $a,a'\in A$ and $b,b'\in B$.
We will write $\Psi(a\otimes b) = \sum_\nu b_\nu\otimes a^\nu$, i.e, $ab
= \sum_\nu b_\nu a^\nu$. All this means that the structure of an algebra
factorisation $X(B,A)$ is fully described by three maps: product in $A$
 $\mu_A: A\otimes
A\to A$, product in $B$ $\mu_B:B\otimes B\to B$ and the twisting
$\Psi: A\otimes B\to B\otimes A$. 

A {\em deformation of an algebra factorisation} $X(B,A)$ over $k$ is an algebra
factorisation $X_t(B_t,A_t)$ over $k[[t]]$ such that the algebras $A_t$,
$B_t$ and $X_t$ are
algebra deformations of $A$, $B$ and $X$ respectively. 
This means that each of the maps $\mu_{A_t}$,
$\mu_{B_t}$ and $\Psi_t$
corresponding to $X_t(B_t,A_t)$ can be written as a formal power series
\begin{equation}
\mu_{A_t} = \mu_A +\sum_{i=1}t^i\mu_A^\su i,\qquad \mu_{B_t} = \mu_B
 +\sum_{i=1}t^i\mu_B^\su
i,\qquad \Psi_t = \Psi +\sum_{i=1}t^i\Psi^\su i,
\label{expan}
\end{equation}
where $\mu_A^\su i :A\otimes A\to A$, $\mu_B^\su i :B\otimes B\to B$,
$\Psi^\su i :A\otimes B\to B\otimes A$, and $\mu_A$, $\mu_B$ and 
$\Psi$ describe
 the factorisation $X(B,A)$. This definition of a deformation of an
algebra factorisation generalises the definition introduced in
\cite{CapSch:twi}, where only the map $\Psi$ was allowed to be deformed.
The need for 
such a generalisation arises from
 the theory of quantum and coalgebra principal bundles 
\cite{BrzMaj:coa}. As
explained in \cite{BrzMaj:fac} the structure of a classical principal
bundle is encoded in the factorisation built on the 
algebra of functions on
the total space of a bundle and the group algebra of the structure
group. The action of the structure group determines the twisting $\Psi$.
 Similarly the structure of a coalgebra 
principal bundle is encoded into an
algebra factorisation, termed a {\em Galois factorisation} in 
\cite{BrzMaj:fac}. In many cases the algebras on which this
factorisation is built are deformations of their classical counterparts
so that not only the twisting map $\Psi$ but also the algebras $A$ and
$B$ are deformed. 

In this note we show that, similarly to the
Gerstenhaber theory of deformation of algebras \cite{Ger:def}, there is a
cohomological interpretation of deformations of algebra factorisations.
Interestingly, such an interpretation uses the
total cohomology of a certain double complex. The situation is therefore
somewhat reminiscent of 
 the Gerstenhaber-Schack theory of deformations of bialgebras
\cite{GerSch:bia}. This is not entirely surprising as there is a close
relationship between algebra factorisations and {\em entwining
structures} (cf. \cite{BrzMaj:coa}). The latter can be seen as a
generalisation of a bialgebra, and, from this point of view, a need for
a double complex in the description of algebra factorisations should be
expected. Deformation theory of entwining structures as well as the
corresponding cohomology theory will be discussed elsewhere.

\section{Construction of the cochain complex}
The fact that $X(B,A)$ is a 
factorisation implies that $B\otimes A^n$ is a right $B$-module via
application of $\Psi$ $n$-times. $A^n$ here denotes the $n$-fold tensor
product of $A$. This in turn implies that $B\otimes
A^n$ is an $(X,X)$-bimodule with the following structure maps. Left
action is obtained by viewing $B\otimes A^n$ as $X\otimes A^{n-1}$ and
multiplying from the left by elements of $X$. The right action is
obtained by viewing $X$ as $B\otimes A$ and acting on $B\otimes
A^n$ by $B$ as described above and then multiplying last factor by
elements in $A$. Similarly (by interchanging $A$ with $B$ and ``left" with
``right"), one makes $B^n\otimes A$ into a left $A$-module
and then an $(X,X)$-bimodule. Using this bimodule structure of $B\otimes
A^n$ and $B^n\otimes A$ one constructs the cohomology of the
factorisation $X(B,A)$ as follows.

First recall that the bar resolution of an  algebra $A$ is a chain
complex $\Bar(A) = ({\rm Bar}_\bullet(A),\delta_A)$, where
$$
{\rm Bar}_n(A) = A^{n+2},\quad \delta_{A,n} =
\sum_{k=0}^{n} (-1)^k A^{k}\otimes\mu_A\otimes A^{n-k} :
A^{n+2}\to A^{n+1}.
$$
Consider bar resolutions of $A$ and $B$. Apply functor $-\otimes A$ to
$\Bar(B)$ and the functor $B\otimes -$ to $\Bar(A)$. Since the
definition of a bar complex boundary operator involves the product in
the algebra only one easily finds that both $\Bar(B)_\bullet\otimes A$,
$B\otimes \Bar(A)_\bullet$ are $(X,X)$-bimodules and 
$\delta_B\otimes A$, $B\otimes\delta_A$ are bimodule maps. This implies that for any
$(X,X)$-bimodule 
$M$ there is a double cochain complex 
$$
C(X(B,A),M) = {}_X{\rm Hom}_X((B\otimes \Bar(A))\otimes _X(\Bar(B)\otimes
A), M).
$$
Explicitly, the space of $(m,n)$-cochains is 
$$
C^{m,n}(X(B,A),M) = {}_X{\rm Hom}_X((B\otimes A^{m+2})\otimes _X
(B^{n+2}\otimes
A), M) \cong {\rm Hom}(A^m\otimes B^n,M).
$$
This last identification is obtained as follows: for
each $\phi\in {}_X{\rm Hom}_X((B\otimes A^{m+2})\otimes _X
(B^{n+2}\otimes
A), M) $ one defines $f_\phi\in {\rm Hom}(A^m\otimes B^n,M)$ via
$$
f_\phi(a_m,\ldots,a_1,b^1,\ldots, b^n) =
\phi(1_B,1_A,a_m,\ldots,a_1,1_A,1_B,
b^1,\ldots, b^n, 1_B,1_A),
$$
while for any $f\in {\rm Hom}(A^m\otimes B^n,M)$ one defines 
$\phi_f\in {}_X{\rm Hom}_X((B\otimes A^{m+2})\otimes _X
(B^{n+2}\otimes
A), M) $ by
\begin{eqnarray*}
&&\phi_f(b,a_{m+1},\ldots,a_0,b^0,\ldots, b^{n+1},a)\\
&& \hspace{.5cm}
=\sum_{(\nu)(\tilde{\nu})}(bb^0_{\nu_0\ldots\nu_{m+1}}a_{m+1}^{\nu_{m+1}})\cdot
 f( a_{m}^{\nu_m},\ldots ,a_1^{\nu_1},b^1_{\tilde{\nu}_1},\ldots ,
b^n_{\tilde{\nu}_n})\cdot
(b^{n+1}_{\tilde{\nu}_{n+1}}a_{0}^{\nu_0\tilde{\nu}_1\ldots
\tilde{\nu}_{n+1}}a).
\end{eqnarray*}
The multi-index notation used above refers to multiple applications of
the twisting map $\Psi$, i.e.,
$$
\sum_{(\nu)}b_{\nu_0\ldots\nu_{n}}\otimes a_{n}^{\nu_{n}}\otimes
\cdots \otimes a_1^{\nu_1} = (\Psi\otimes A^{n-1})(A\otimes \Psi\otimes
A^{n-2})\cdots (A^{n-1}\otimes\Psi)(a_{n}\otimes\cdots\otimes
a_1\otimes b),
$$
etc. From the deformation theory point of view the case $M=X$ is of the
greatest interest, thus to this case we restrict our attention from now 
on. Consider the
Hochschild complexes of $A$ and $B$. Notice that the  inclusions
${\rm Hom}(A^m,A)\hookrightarrow {\rm Hom}(A^m,X)$ and ${\rm Hom}(B^n,B)
\hookrightarrow {\rm Hom}(B^n,X)$ given by $f\mapsto 1_B\otimes f$ and
$g\mapsto g\otimes 1_A$ respectively are inclusions of cochain 
complexes. The complex obtained from $C(X(B,A),X)$ by replacing
$C^{\bullet, 0}(X(B,A),X)$ by the Hochschild complex of $A$, 
and $C^{0,\bullet}(X(B,A),X)$ by the Hochschild complex of $B$ 
 is denoted by $C(X(B,A))$. Explicitly, one has the following double
complex 

{\tiny
$$
\begin{CD}
B\oplus A @>{d_A}>> {\rm Hom}(A,A) @>{d_A}>> {\rm Hom}(A^2,A) @>{d_A}>> 
{\rm Hom}(A^3,A) @>{d_A}>>\\
@VV{d_B}V  @VV{d_B}V    @VV{d_B}V   @VV{d_B}V\\
{\rm Hom}(B,B) @>{d_A}>> {\rm Hom}(A\otimes B,X) @>{d_A}>> {\rm Hom}
(A^2\otimes B, X) @>{d_A}>> {\rm Hom}
(A^3\otimes B, X) @>{d_A}>>\\
@VV{d_B}V   @VV{d_B}V   @VV{d_B}V   @VV{d_B}V\\
{\rm Hom}(B^2, B) @>{d_A}>> {\rm Hom}(A\otimes B^2, X) @>{d_A}>> {\rm
Hom}(A^2\otimes B^2, X) @>{d_A}>>{\rm Hom}
(A^3\otimes B^2, X) @>{d_A}>>\\
@VV{d_B}V   @VV{d_B}V   @VV{d_B}V   @VV{d_B}V
\end{CD}
$$}
The coboundary operators $d_A:C^{m,n}(X(B,A))\to C^{m+1,n}(X(B,A))$ and 
$d_B:C^{m,n}(X(B,A))\to C^{m,n+1}(X(B,A))$, $m,n>0$
are given explicitly by 
\begin{eqnarray}
d_Af(a_{m+1},\ldots,a_1,b^1,\ldots, b^n) &=& 
a_{m+1}f(a_{m},\ldots,a_1,b^1,\ldots,b^n) \nonumber \\
&& \hspace{-1.5cm} + \sum_{i=0}^{m-1}(-1)^{i+1}f(a_{m+1},\ldots,a_{m+1-i}a_{m-i},
\ldots ,a_1,b^1,\ldots, b^n)\nonumber \\
&& \hspace{-1.5cm} + (-1)^{m+1}\sum_{(\nu)}f(a_{m+1},\ldots,a_2,b^1_{\nu_1},\ldots,
b^n_{\nu_n}) a_1^{\nu_1\ldots \nu_{n}},
\label{da}
\end{eqnarray}
and 
\begin{eqnarray}
d_Bf(a_{m},\ldots,a_1,b^1,\ldots, b^{n+1}) &=& 
\sum_{(\nu)} b^{1}_{\nu_1\ldots \nu_{m}}
f(a_{m}^{\nu_m},\ldots,a_1^{\nu_1},b^2,\ldots,b^{n+1}) \nonumber\\
&& \hspace{-1.5cm} +
\sum_{i=1}^{n}(-1)^{i}f(a_{m},\ldots,a_1,b^1,\ldots,
b^{i}b^{i+1},\ldots, b^{n+1})\nonumber \\
&& \hspace{-1.5cm} + (-1)^{n+1}f(a_{m},\ldots,a_1,b^1,\ldots,
b^n) b^{n+1}.
\label{db}
\end{eqnarray}
For $n=0$, $d_A$ is  the usual Hochschild coboundary, while the 
$d_B$ are given by (\ref{db}), provided one views ${\rm
Hom}(A^m,A)$ inside ${\rm
Hom}(A^m,X)$ first. Similarly, for 
$m=0$, $d_B$ is  the usual Hochschild coboundary, while the 
$d_A$ are given by (\ref{da}), provided one views ${\rm
Hom}(B^n,B)$ inside ${\rm
Hom}(B^n,X)$ first.
The construction of the above complex implies immediately that $d_A\circ
d_B = d_B\circ d_A$, so that one can combine the double
complex into a complex
$(C^\bullet(X(B,A)),D)$,
$$
C^n(X(B,A)) = {\rm Hom} (A^n,A)\oplus \bigoplus_{k=1}^{n-1} {\rm
Hom}(A^{n-k}\otimes B^k, X)\oplus {\rm Hom} (B^n,B),
$$ 
$ D\mid_{C^{m,n}} = (-1)^md_B +d_A.
$
The cohomology of the complex $(C^\bullet(X(B,A)),D)$ is denoted
by $H^\bullet(X(B,A))$. 

\section{Cohomological interpretation of deformations}
Two deformations $X_t(B_t,A_t)$ and $\tilde{X}_t(\tilde{B}_t,\tilde{A}_t)$ of an algebra
factorisation $X(B,A)$ are said to be {\em equivalent} to each other
 if there exist algebra isomorphisms
$\alpha_t :A_t\to \tilde{A}_t$, $\beta_t:B_t\to \tilde{B}_t$ 
of the form $\alpha_t = A +
\sum_{i=1} t^i\alpha^\su i$, $\beta_t = B +
\sum_{i=1} t^i\beta^\su i$, and such that $\beta_t\otimes \alpha_t:
X_t\to\tilde{X}_t$ is an algebra isomorphism.  A deformation
$X_t(B_t,A_t)$ is called a {\em trivial deformation} if it is 
equivalent to an
algebra factorisation in which all the maps $\mu_A^\su i$, $\mu_B^\su
i$, $\Psi^\su i$ in (\ref{expan}) vanish. 
An {\em infinitesimal deformation} of $X(B,A)$ is a
deformation of $X(B,A)$ modulo $t^2$. 
\begin{theorem}
 There is a one-to-one correspondence between the equivalence classes of
infinitesimal deformations of $X(B,A)$ 
and  $H^2(X(B,A))$.
\end{theorem}
\begin{proof}
 For an infinitesimal deformation it is enough  to consider $\mu_{A_t}
= \mu_A+ t\mu_A^\su 1$, $\mu_{B_t}
= \mu_B+ t\mu_B^\su 1$, $\Psi_t
= \Psi+ t\Psi^\su 1$  where $\mu_A^\su 1 \oplus \Psi^\su 1 \oplus
\mu_B^\su 1 \in C^{2,0}\oplus C^{1,1}\oplus C^{0,2} = C^{2}(X(B,A))$. First we show that the triple
$
(\mu_A^\su 1,\Psi^\su 1,
\mu_B^\su 1 )$ defines an infinitesimal deformation if and only if 
$\mu_A^\su 1 \oplus \Psi^\su 1 \oplus
\mu_B^\su 1 $ is a cocycle.

As the first row and the first column in $C(X(B,A))$ are simply
 Hochschild complexes, 
a standard algebra deformation theory argument shows that the
associativity of $\mu_{A_t}$ and $\mu_{B_t}$ modulo $t^2$ is equivalent to
the conditions $d_A\mu_A^\su 1 = d_B\mu_B^\su 1 = 0$. 
In view of this fact we need to
show that $\Psi_t$ satisfies conditions (\ref{diag.A}) and
(\ref{diag.B}) if and only if $d_B\mu_A^\su 1 + d_A\Psi^\su 1 =0$  and 
$d_A\mu_B^\su 1 - d_B\Psi^\su 1 =0$. Expanding
(\ref{diag.A}), (\ref{diag.B}) for $\Psi_t$ in powers of $t$ one easily 
finds that the $t^0$-order terms are simply equations (\ref{diag.A}),
(\ref{diag.B}) for $\Psi$. Therefore only terms of order $t$ require
further study.  The
$t$-order term in the first of equations (\ref{diag.A}) is
\begin{eqnarray*}
&& \hspace{-31pt}(B\otimes \mu_A)(\Psi^\su 1\otimes A)(A\otimes \Psi) - 
\Psi^\su 1(\mu_A\otimes B) + (B\otimes\mu_A)(\Psi\otimes A)
(A\otimes \Psi^\su 1)\\
&& - \Psi(\mu_A^\su 1\otimes B) + (B\otimes\mu_A^\su 1)(\Psi\otimes A)
(A\otimes \Psi) =0 .
\end{eqnarray*}
This is precisely the statement that $d_B\mu_A^\su 1 +d_A\Psi^\su 1 = 0$. 
Evaluating this condition at
$1_A\otimes 1_A\otimes b$ one easily finds that 
$\Psi^\su 1 (1_A\otimes b) =
-b\otimes \mu_A^\su 1(1_A\otimes 1_A) +\Psi(\mu_A^\su 1(1_A\otimes 
1_A)\otimes b)$, 
i.e., $\Psi_t(1_{A_t}\otimes b) = b\otimes 1_{A_t}$, where $1_{A_t} =
1_A - \mu_A^\su 1(1_A\otimes 1_A)t$ is the unit in the infinitesimal
deformation $A_t$. Thus the second of
equations (\ref{diag.A}) holds. Similarly one shows that the equations
(\ref{diag.B}) hold for $\Psi_t$ modulo $t^2$ if and only if
$d_A\mu_B^\su 1-d_B\Psi^\su 1 = 0$. Therefore the necessary and sufficient condition
for $X(B,A)_t$ to be an infinitesimal deformation of
$X(B,A)$ is that $\mu_A^\su 1 \oplus \Psi^\su 1\oplus \mu_B^\su 1$ 
be a 2-cocycle in $C^2(X(B,A))$ as
required.

Let $X_t(B_t,A_t)$ and $\tilde{X}_t(\tilde{B}_t,\tilde{A}_t)$ be two
infinitesimal deformations of an algebra
factorisation $X(B,A)$ given by the cocycles $\mu_A^\su 1 \oplus \Psi^\su 1\oplus \mu_B^\su
1$ and $\tilde{\mu}_A^\su 1 \oplus \tilde{\Psi}^\su 1\oplus 
\tilde{\mu}_B^\su 1$ respectively. We need to show that these two
 deformations
are equivalent to each other modulo $t^2$ if and only if the corresponding 
cocycles differ by a coboundary.
In view of the Gerstenhaber theory,  $\alpha_t =
A+t\alpha :A_t\to \tilde{A}_t$ and  $\beta_t =
B+t\beta :B_t\to \tilde{B}_t$ are the algebra isomorphisms modulo $t^2$
if and only if  ${\mu}_A^\su 1-\tilde{\mu}_A^\su 1 = d_A\alpha$ and
${\mu}_B^\su 1-\tilde{\mu}_B^\su 1 = d_B\beta$. Thus it remains to be shown that
$\beta_t\otimes\alpha_t :X_t\to \tilde{X}_t$ is an algebra isomorphism
modulo $t^2$ if and only if 
$$
{\Psi}^\su 1 - \tilde{\Psi}^\su 1= d_A\beta - d_B\alpha.
$$

Suppose that  $\phi_t = \beta_t\otimes\alpha_t$ is an 
isomorphism of algebras,
and let $\phi$ be the $t$-order term in the expansion of $\phi_t$, i.e.,
\begin{equation}
\phi = B\otimes \alpha
+\beta\otimes A.
\label{phi}
\end{equation}
Since 
$\phi_t$ is an algebra map
\begin{equation}
\phi_t(ab) = \phi_t(a)\phi_t(b), \qquad \forall a\in A, b\in B.
\label{alg}
\end{equation}
Note that the product on the left hand side of (\ref{alg}) is in
$X_t(B_t,A_t)$ while on the right hand side is in 
$\tilde{X}_t(\tilde{B}_t,\tilde{A}_t)$. One
can use (\ref{phi}) to find the 
$t$-order term on the left hand side of (\ref{alg}) 
$$
\sum_\nu\phi(b_\nu a^\nu) +\Psi^\su 1(a\otimes b) = \sum_\nu b_\nu
\alpha(a^\nu) +\sum_\nu\beta(b_\nu)a^\nu + \Psi^\su 1(a\otimes b).
$$
All the products are in $X(B,A)$ now. 
On the other hand the $t$-order term on the right hand side is
$$
\tilde{\Psi}^\su 1(a\otimes b) +a\phi(b) +\phi(a)b = 
\tilde{\Psi}^\su 1(a\otimes b) +a\beta(b)
+\alpha(a)b.
$$
Thus if $\phi_t$ is an algebra map we have
\begin{eqnarray*}
({\Psi}^\su 1 - \tilde{\Psi}^\su 1)(a\otimes b) & = & 
a\beta(b)-\sum_\nu\beta(b_\nu)a^\nu
 +\alpha(a)b -\sum_\nu b_\nu
\alpha(a^\nu) \\
&=& (d_A\beta - d_B\alpha)(a\otimes b),
\end{eqnarray*}
 as required.

To prove the converse one needs to repeat the same computations in
reversed order.
\end{proof}

The next step usually undertaken in the deformation theory, 
is to study obstructions for extending a deformation
modulo $t^{n}$ to a deformation modulo $t^{n+1}$. Such an obstruction
consists of four terms. The first two terms come from the deformation
of algebra structures of $A$ and $B$, one for each algebra. They are:
$$
\ob_A^\su n = \sum_{k=1}^{n-1}[\mu_A^\su k(\mu_A^\su{n-k}\otimes A) - \mu_A^\su
k(A\otimes \mu_A^\su{n-k})], 
$$
$$ 
\ob_B^\su n = \sum_{k=1}^{n-1}[\mu_B^\su k(\mu_B^\su{n-k}\otimes B) - \mu_B^\su
k(B\otimes \mu_B^\su{n-k})], 
$$
The
remaining two obstructions arise from the factorisation conditions
(\ref{diag.A}) and (\ref{diag.B}):
$$
\ob_{A,\Psi}^\su n = 
\sum_{k=1}^{n-1}\Psi^\su{n-k}(\mu_A^\su k\otimes B) -  
\sum_{\stackrel{\scriptstyle k,l=0}{1\leq k+l\leq n}}^{n-1}
(B\otimes \mu_A^\su k)(\Psi^\su
l\otimes A)(A\otimes \Psi^\su{n-k-l}),
$$
$$
\ob_{B,\Psi}^\su n = 
\sum_{\stackrel{\scriptstyle k,l=0}{1\leq k+l\leq n}}^{n-1}
  (\mu_B^\su k\otimes A)
(B\otimes\Psi^\su l)(\Psi^\su{n-k-l}\otimes B) -
\sum_{k=1}^{n-1}\Psi^\su{n-k}(A\otimes\mu_B^\su k).
$$
Here $\Psi^\su 0 = \Psi$, $\mu_A^\su 0 = \mu_A$ and $\mu_B^\su 0 =
\mu_B$. The following theorem is an algebra factorisation version of a 
standard result in the deformation
theory. 
\begin{theorem}
If $X_t(B_t,A_t)$ is a deformation of $X(B,A)$ modulo $t^{n}$ then
$$
\ob^\su n = \ob_A^\su n\oplus\ob_{A,\Psi}^\su 
n\oplus\ob_{B,\Psi}^\su n\oplus\ob_B^\su n 
$$
is a 3-cocycle in the complex $C(X(B,A))$. 
$X_t(B_t,A_t)$ can be extended to 
a deformation of $X(B,A)$ modulo $t^{n+1}$ if and only if $\ob^\su n$ is
a coboundary.
\label{thm.obs}
\end{theorem}
\begin{proof}
The first part of the theorem can be proven in the following way 
(standard in the
deformation theory of algebras, which also asserts  that 
$\ob_A^\su n$ and
$\ob_B^\su n$ are Hochschild cocycles).  Let 
$$
\tilde{\mu}_A= \mu_A +\sum_{i=1}^{n-1}t^i\mu_A^\su i,\qquad
\tilde{\mu}_B = \mu_B
 +\sum_{i=1}^{n-1}t^i\mu_B^\su
i,\qquad \tilde{\Psi} = \Psi +\sum_{i=1}^{n-1}t^i\Psi^\su i.
$$
The proof hinges on two observations. Firstly, one easily finds that
$$
\ob_A^\su n = \mbox{\rm coefficient of $t^n$ in}\;  \tilde{\mu}_A(\tilde{\mu}_A\otimes
A) - \tilde{\mu}_A(A\otimes \tilde{\mu}_A),
$$
$$
\ob_{A,\Psi}^\su n = \mbox{\rm coefficient of $t^n$ in}\; \tilde{\Psi}
(\tilde{\mu}_A\tens B)-(B\tens\tilde{\mu}_A)(\tilde{\Psi} \otimes
A)(A\otimes \tilde{\Psi}),
$$
$$
\ob_{B,\Psi}^\su n = \mbox{\rm coefficient of $t^n$ in} \;
(\tilde{\mu}_A\tens B)(B\otimes
\tilde{\Psi})(\tilde{\Psi}\otimes B) - \tilde{\Psi}(A\tens\tilde{\mu}_B),
$$
$$
\ob_B^\su n = \mbox{\rm coefficient of $t^n$ in }\; \tilde{\mu}_B(\tilde{\mu}_B\otimes
B) - \tilde{\mu}_B(B\otimes \tilde{\mu}_B).
$$
Secondly one should notice that 
$$
D\ob^\su n = \mbox{\rm coefficient of $t^n$ in }\; \tilde{D} \ob^\su n,
$$
where $\tilde{D}$ is obtained by replacing $\mu_A$, $\mu_B$ and $\Psi$
in definition of $D$ with  $\tilde{\mu}_A$, $\tilde{\mu}_B$ and 
$\tilde{\Psi}$. Expanding $\tilde{D} \ob^\su n$, with $\ob^\su n$
expressed entirely in terms of the tilded structure maps, 
one discovers that the 
term-by-term cancellations yield $\tilde{D} \ob^\su n=0$. 
Thus $\ob^\su n$ is a cocycle as asserted. (This expansion is a
straightforward procedure, one only has to remember to take the 
inclusions of Hochschild cocycles into $C(X(B,A))$ properly into
account.) 
 
It follows from the Gerstenhaber theory that
 $A_t$ and $B_t$ are deformations of $A$ and
$B$ respectively modulo $t^{n+1}$ if and only if $\ob_A^\su n$ and
$\ob_B^\su n$ are coboundaries in the Hochschild cohomology, i.e.,
there exist $\mu_A^\su n: A\otimes A\to A$ and $\mu_B^\su n:B\otimes
B\to B$ such that $d_A\mu_A^\su n = \ob_A^\su n$ and $d_B\mu_B^\su n = \ob_B^\su n$. 
 Thus only
the conditions arising from (\ref{diag.A}) and (\ref{diag.B}) require
further study.
 Gathering all the terms of order $t^n$ in (\ref{diag.A}) and
(\ref{diag.B}) one easily finds  that $X_t(B_t,A_t)$ is a deformation
modulo $t^{n+1}$ if and only if
$$
d_A\Psi^\su n  +  d_B\mu_A^\su n = \ob_{A,\Psi}^\su n, \qquad 
d_A\mu_B^\su n - d_B\Psi^\su n = \ob_{B,\Psi}^\su n.
$$
All this means that the necessary and sufficient condition for 
$X_t(B_t,A_t)$ to be a deformation
modulo $t^{n+1}$ is that
$$
D(\mu_A^\su n\oplus\Psi^\su n\oplus \mu_B^\su n) =
\ob^\su n,
$$
i.e., $\ob^\su n$ is a coboundary, as required.
\end{proof}

An interesting special case of this general deformation theory is 
a deformation $X_t(B,A)$, i.e., the
algebras $B$ and $A$ are not deformed, and only the formal power series 
$\Psi_t$ is non-trivial. This type of deformation is
considered in \cite{CapSch:twi}. In this case $\ob_A^\su n = 0$, $\ob_B^\su n = 0$, and
$$
\ob_{A,\Psi}^\su n =  -  
\sum_{i=1}^{n-1}
(B\otimes \mu_A)(\Psi^\su
i\otimes A)(A\otimes \Psi^\su{n-i}),
$$
$$
\ob_{B,\Psi}^\su n = 
\sum_{i=1}^{n-1}
  (\mu_B\otimes A)
(B\otimes\Psi^\su i)(\Psi^\su{n-i}\otimes B).
$$
The obstruction removing equations are:
$$
d_A\Psi^\su n  = \ob_{A,\Psi}^\su n, \qquad 
d_B\Psi^\su n = -\ob_{B,\Psi}^\su n,
$$
and coincide with the equations given in
\cite[Theorem~4.11]{CapSch:twi}. 

We would like to conclude with three concrete examples
illustrating the deformation theory of algebra factorisations. The first
example deals with a deformation affecting both the algebra structure of
$A$ as well as the map $\Psi$, while the remaining two are an
illustration of a deformation of $\Psi$ only.

\begin{example}
\rm 
 Let $k =\bf C$, $A= {\bf C}[a,\ba ]/(a\ba-\ba a)$ 
(an algebra of polynomials in two commuting variables), $B={\bf C}[b]$, 
and 
$X= B\otimes A$ a tensor product algebra (an algebra of polynomials in
three commuting variables). $A$ is spanned by all 
monomials $a^k\ba^l$, $k,l=0,1,2,\ldots$, while 
$B$ is spanned by 
the set $\{ b^r\; |\; r=0,1,2,\ldots \}$. Thus the structure maps 
are $\mu_A(a^k\ba^l\otimes a^r\ba^s) = a^{k+r}\ba^{l+s}$, 
$\mu_B(b^r\otimes b^s)
=b^{r+s}$, while the map $\Psi$ is the usual
twist, $\Psi(a^k\ba^l\otimes b^r) = b^r\otimes a^k\ba^l$. 
One easily verifies
that $\mu_A^\su 1\oplus\Psi^\su 1$, where
$$
\Psi^\su 1(a^k\ba^l\otimes b^r) = lrb^r\otimes a^k\ba^l - 
krb^{r+1}\otimes a^{k-1}\ba^{l+1},
\quad \mu_A^\su 1(a^k\ba^l\otimes a^r\ba^s) = lr a^{k+r}\ba^{l+s},
$$
is a cocycle in $C^\bullet(X(B,A))$ and therefore defines an infinitesimal
deformation of $X(B,A)$.  The $n=2$ obstruction 3-cocycle consists of three
terms:
\begin{eqnarray*}
&&\ob_A^\su 2 (a^k\ba^l\otimes a^m\ba^n\otimes a^p\ba^r) =
lp(lm-np)a^{k+m+p}\ba^{l+n+r},\\
&&\ob_{A,\Psi}^\su 2(a^k\ba^l\otimes a^m\ba^n\otimes b^r) = 
r((lm+kn)r+km)b^{r+1}\otimes
a^{k+m-1}\ba^{l+n+1}\\
&&\hspace{1.2in} -lnr^2b^r\otimes a^{k+m}\ba^{l+n}
-kmr(r+1)b^{r+2}\otimes
a^{k+m-2}\ba^{l+n+2},\\
&&\ob_{B,\Psi}^\su 2(a^k\ba^l\otimes b^r\otimes b^s) = l^2rsb^{r+s}\otimes
a^k\ba^l -(2l+1)krsb^{r+s+1}\otimes
a^{k-1}\ba^{l+1}\\
&&\hspace{1.2in} + (k-1)krs b^{r+s+2}\otimes
a^{k-2}\ba^{l+2}.
\end{eqnarray*}
This obstruction can be removed by setting
\begin{eqnarray*}
&&\mu_A^\su 2(a^k\ba^l\otimes a^m\ba^n) = lm(c+{lm\over
2})a^{k+m}\ba^{l+n} \\
&&\Phi^\su 2(a^k\ba^l\otimes b^r) = lr(c+{lr\over 2})b^r\otimes
a^k\ba^l - kr({k+r-1\over 2} +lr +c)b^{r+1}\otimes
a^{k-1}\ba^{l+1}\\
&&\hspace{1.2in} {1\over 2}k(k-1)r(r+1)b^{r+2}\otimes
a^{k-2}\ba^{l+2},
\end{eqnarray*}
where $c$ is a constant number. This deformation of $X(B,A)$ modulo 
$t^3$ has 
presentation with generators $a$, $\ba$ and $b$, and the relations
$$
\ba a = qa\ba, \quad \ba b=qb\ba, \quad ab=ba+(1-q)b^2\ba,
$$
where $q=1+t+(c+{1\over 2})t^2$. 
A simple calculation reveals that the above relations define an
associative algebra to all powers of $t$. Furthermore the same
calculation implies that the above relations with $q=1+t+(c+{1\over 2})t^2+ o(t^3)$ 
 describe an associative algebra factorisation $X_t(B,A_t)$ over 
${\bf C}[[t]]$. The deformed product in $A$ and deformed twisting map 
come out as:
$$
\mu_{A_t}(a^k\ba^l\otimes a^r\ba^s) = q^{lr} a^{k+r}\ba^{l+s},
$$
$$
\Psi_t(a^k\ba^l\otimes b^r) = q^{lr}\sum_{i=0}^k\binom{k}{i}_q
\binom{r+i-1}{i}_q(q;q)_i b^{r+i}\otimes a^{k-i}\ba^{l+i},
$$
where
$(q;q)_{-1}=(q;q)_0 =1$, $(q;q)_i = (1-q)(1-q^2)\cdots (1-q^i)$, 
$$
\binom{k}{i}_q= \frac{(q;q)_k}{(q;q)_i(q;q)_{k-i}}.
$$
\label{ex.first}
\end{example}

\begin{example}
\rm Let $k=\bf C$, $A = {\bf C}[a,\ba]/(\ba a-qa\ba)$ (Manin's
quantum plane), $B= {\bf C}[b]$ and $\Psi:A\otimes B\to B\otimes A$, 
$\Psi(a^k\ba^l\otimes b^r) = q^{lr}b^r\otimes a^k\ba^l$, 
where $q$ is a non-zero complex
number which is not a root of unity. The resulting algebra
factorisation is
$$
X(B,A) = {\bf C}[b,a,\ba]/(\ba a-qa\ba, \ba b  - qb\ba, ab-ba).
$$
 One easily verifies
that $\Psi^\su 1:A\otimes B\to B\otimes A$,
$$
\Psi^\su 1(a^k\ba^l\otimes b^r) = q^{lr}[k][r]
b^{r+1}\otimes a^{k-1}\ba^{l+1},
$$
where $[k] = \frac{1-q^k}{1-q}$, 
is a cocycle in $C^\bullet(X(B,A))$ and therefore defines an infinitesimal
deformation of $X(B,A)$. This infinitesimal deformation of $X(B,A)$ 
has a presentation with generators $a,\ba,b$ and relations:
\begin{equation}
\ba a = qa\ba, \quad \ba b=qb\ba, \quad ab=ba+tb^2\ba.
\label{ex2.eq}
\end{equation}
The infinitesimal deformation given by relations (\ref{ex2.eq}) 
can be extended to a
deformation to all orders in $t$ by setting:
$$
\Psi^\su i(a^k\ba^l\otimes b^r) = q^{lr}[k][k-1]\cdots[k-i+1]
\binom{r+i-1}{i}_q b^{r+i}\otimes a^{k-i}\ba^{l+i},
$$
with the same notation as in Example~\ref{ex.first} (but note that $q$
has a different meaning). The resulting algebra factorisation 
$X_t(B,A)$ over ${\bf C}[[t]]$ has presentation with the relations
(\ref{ex2.eq}).  The fact that with these definitions of the $\Psi^\su i$
 all the
obstructions are removed can be verified directly by using various
identities for $q$-binomial coefficients. This can also be verified by
checking directly that relations (\ref{ex2.eq}) define an associative 
algebra. 
\label{ex.second}
\end{example}

\begin{example}
Deformation of quaternions.
\rm  Let $k= {\bf R}$, 
$A = {\bf R}[i]/(i^2 +1) \cong {\bf C}$, $B = {\bf R}[j]/(j^2 +1) \cong {\bf C}$, 
and set $\Psi(i\otimes j) = -j\otimes i$. The resulting algebra is
$X = {\bf H}$, i.e., there is a factorisation ${\bf H}({\bf C},{\bf
C})$. It is well known that the second Hochschild cohomology of ${\bf
C}$ viewed as a real algebra is trivial. Thus we can choose $\mu_A^\su 1
=0$ and $\mu_B^\su 1 =0$. This implies that any cocycle in $C^2({\bf H}({\bf C},{\bf
C}))$ must be cohomology equivalent to the map $\Psi^\su 1:{\bf C}\otimes {\bf
C}\to {\bf C}\otimes {\bf C}$ such that $\Psi^\su 1(1\otimes j)
= \Psi^\su 1(i\otimes 1)
= \Psi^\su 1(1\otimes 1)=0$. All such cocycles  can be easily computed, and
one finds that  ${\rm dim}H^2({\bf H}({\bf C},{\bf C})) =
1$ with the unique cohomology class generated by the cocycle 
$\Psi^\su 1(i\otimes j) = 1\otimes 1$. Consequently, 
 there is only one deformation of quaternions
that retains the factorisation property, namely the algebra factorisation
${\bf H}_t({\bf C},{\bf C})$ 
spanned by $1$, $i$, $j$ and $ij$ subject to the relations 
$i^2 = j^2 =-1$, $ij+ji = t$. This result agrees with the classification
of factorisations given in 
\cite[Example~2.12]{CaeIon:fac}. 
\end{example}
\begin{center}
ACKNOWLEDGMENTS
\end{center}
This paper was written when the author was a Lloyd's of London 
Tercentenary
Research Fellow at the Department of Mathematics, University of York.

\end{document}